\documentclass[12pt,a4paper]{article}
\usepackage{color,amsmath}
\usepackage{amsfonts}
\usepackage{amssymb}
\usepackage{graphicx}

\newtheorem{theorem}{Theorem}

\begin{document}

\

\begin{center}
{\Large \textbf{Selection of variables for cluster analysis and classification
rules}}

\

{\large  Ricardo Fraiman$^{\ast}$, Ana Justel$^{\ast\ast}$\footnote{Corresponding
author: Ana Justel, Departamento de Matem\'aticas, Universidad Aut\'onoma de Madrid.
Campus de Cantoblanco, 28049 Madrid, Spain. Email: ana.justel@uam.es} and Marcela
Svarc$^{\ast}$}

\

\noindent\emph{$^{\ast}$Departamento de Matem\'atica y Ciencias, Universidad de San
Andr\'es, Argentina} \\[0pt]\emph{$^{\ast\ast}$Departamento de Matem\'aticas,
Universidad Aut\'onoma de Madrid, Spain}
\end{center}

\

\begin{abstract}
In this paper we introduce  two procedures for variable selection in cluster
analysis and classification rules. One is mainly oriented to detect the ``noisy''
non--informative variables, while the other deals also with multicolinearity. A
forward--backward algorithm is also proposed to make feasible these procedures in
large data sets. A small simulation is performed and some real data examples are
analyzed.
\newline \noindent{\em A.M.S. 1980 subject classification:} {\em Primary: } 62H35
\end{abstract}

\section{Introduction}

In multivariate analysis there are several statistical procedures whose output is a
partition of the space. Typical examples of this situation are cluster analysis and
classification rules. In cluster analysis (or un--supervised classification) we look
for a partition of the space into homogeneous groups or clusters (with small
dispersion within groups), that help us to understand the structure of the data.
Several cluster methods have been proposed, such as hierarchical clustering
(Hartigan, 1975), k-means (MacQueen, 1967), k-mediods (Kaufman and Rousseeuw, 1987),
kurtosis based clustering  (Pe\~na and Prieto, 2001). From most of them we get a
partition of the space in disjoint subsets.

Pattern recognition or classification  is about guessing or predicting the unknown
nature of an observation, a discrete quantity such as black or white, one or zero,
sick or healthy. An observation is a collection of numerical measurements such as an
image (which is a sequence of bits, one per pixel), a vector of weather data, or an
electrocardiogram. In classification rules, we have in addition a training sample
for each group, from which we know together with the observation of the random
vector of variables, a label that indicates to which subpopulation it belongs. Then
a {\it classifier} is any map that represents for each new data our guess of the
class, given its associated vector. The map produces a classification rule, that is
also a partition of the space. According to which subset of the partition a new data
belongs, is classified in that class. There is also an extensive literature on
classification rules, such as Fisher's linear discrimination (Fisher, 1936), nearest
neighbor rules (Fix and Hodges, 1951), regression trees-CART (Breiman et al., 1984),
or reduced kernel discriminant analysis (Hern\'andez and Velilla, 2005).

A general problem in cluster or classification is to find structures in a high
dimensional variable space but with small data sets. It is common that in many
practical cases, the amount of variables (that should not be confused with the
amount of information) is too high. This may be due to the presence of several
``noisy'' non--informative variables, and/or redundant information from strongly
correlated variables that may produce multicolinearity. Then the information
contained in the data set could be extracted from a reduced subset of the original
variables.

A difficult task is to find out which variables are ``important'',
where the concept of ``important''  should be related to the
statistical procedure we are dealing with. If we are interested in
cluster analysis, we would like to find the variables that explains
the groups we have found. In this way, a (small) subset of variables
should ``explain'' as best as possible the statistical procedure in
the original space (the high dimensional space). Dimension reduction
techniques (like principal component analysis) will produce linear
combinations of the variables which are difficult to interpret
unless most of the coefficients of the linear combination are
negligent. The variable selection method of Fowlkes, et al. (1988)
shifts the problem to a reduced variable space and looks for new
clusters with less variables. Tandesse, et al. (2005) propose a
Bayesian approach for simultaneously selecting variables and
identifying cluster structures without knowing the number of
clusters. The Bayesian model with latent variables is very useful in
cluster analysis since it produces the most complete output: number
of clusters, data allocation and informative variables. To solve the
model it is necessary to use MCMC methods, in particular
Metropolis-Hastings with Reversible-Jump (Green, 1995), that
introduce an important complexity to the users that are not familiar
with computer programming.

In this paper we propose consistent statistical methods for variable
selection that are easy to use. The variables that explain better
the procedure on the original space help us to understand better the
cluster output, and as a by-product, we find a dimension reduction
procedure that can be used in a new data set for the same problem.
We consider two different proposals based on the idea of
``blinding'' unnecessary variables. To cancel the effect of one
variable, we substitute all the values of that variable by it's
marginal mean in the first proposal and by the conditional mean in
the second proposal. The  marginal mean approach is mainly oriented
to detect the ``noisy'' non--informative variables, while the
conditional mean approach is more related to deal also with
multicolinearity. The first one is simpler and does not require
large sample size as the second one. In practice, we will also need
an algorithm to solve the optimization problem.

In Section 2 we define in precise terms what we understand for a subset of variables
that explains a multivariate partition procedure. Next we define our objective
function and  provide a strongly consistent estimate of the optimal subset. A small
simulation study is also performed. In Section 3 we introduce the proposal based on
the conditional mean  and show the performance in a simulated data set. In Section 4
we describe a forward--backward selection  algorithm that looks for the minimum
subset that explains a fixed percentage of the data assignation to the clusters.
Section 5 is devoted to the analysis of two real data examples  with medium and
large dimensional variable spaces. Section 6 includes some final remarks and the
proofs are given in the Appendix.

\section{Dropping out noise non--informative variables}

Let $X = (X_1,\ldots, X_p)$ be a random vector with distribution $P$. We consider
any statistical procedure whose output is a partition of the space $R^p$. For
instance, this is the case of the population target for most clustering methods or
classification rules. To fix ideas we will concentrate in cluster methods. For a fix
number of clusters $K$, we have a function
$$
f: R^p \rightarrow \{1, \ldots, K\}
$$
which determines to which cluster each single point belongs. We denote the space
partition by $G_k = f^{-1}(k), \hspace{2mm} k=1, \ldots, K$, that satisfies
$$
P\left( \bigcup_{k=1}^K G_k \right) = 1.
$$
For instance, if we consider $k-means$ (with K=2), and $c_1, c_2 \in R^p$ are the
cluster centers, i.e. the values that minimize
$$
E \left( \min(||X-c_1||^2, ||X-c_2||^2) \right),
$$
the set $G_1$ is given by $G_1 = \{x \in R^p: ||x-c_1|| \leq ||x-c_2|| \}$, while
$G_2 = G_1^c$.

If $p$ is large, typically some of the components of vector $X$ are strongly
correlated or might be almost irrelevant for the cluster procedure. Then, if the
information from the noisy variables is removed from our data, we should expect that
their cluster allocations does not change. These means that the data are kept in the
same group as in the original partition. The key point is to notice that the
partition is defined in the original $p$-dimensional space and the input data
requires information from all the variables, included the noisy ones. We propose to
look for the subset of indices $I\subset \{1,\ldots,p\}$ for which the original
partition rule applied to a new ``less informative`` vector $Y^I \in R^p$ built up
from $X$, behaves as close as possible to the procedure when is applied to the
``full information'' vector $X$. The vector $Y^I$ contains the variables from $X$
that are index by $I$, and the rest of the variables with index outside the set $I$
are ``blinded''. A noisy variable means that its probability distribution is almost
the same at all the clusters. This suggests to substitute the information in the
``blinded '' variables by their mean value.

It will depend on the problem (the distribution $P$ of $X$) and on how many
variables $d < p$ we select, the percentage of cluster allocations explained by
them. In practice, we can choose $d$ in order to explain at least a fixed
percentage, for instance, $90\%$, $95\%$ or $100\%$ of the data.

\subsection{Population and empirical objective functions}

We now put our purpose in a precise setup. Given a subset of indices
$$
I = \{i_1, \ldots, i_d\} \subset \{1,\ldots, p\},
$$
we define the vector $Y^I:= Y = (Y_1, \ldots, Y_p)$, where $Y_i = X_i$ if $ i \in I$
and $Y_i = E(X_i)$ otherwise. Note that instead of the expectation $E(X_i)$ we can
use the median of $X_i$ or any other location parameter for the $i$--coordinate like
M-estimates or trimmed--means. The results will still holds provided we have a
strong consistent estimate of the location parameter.

For a fixed integer $d < p$, the population target is the set $I \subset \{1,
\ldots, p \}$, $\#I$ = d, for which the {\em population objective function}, given
by
$$
h(I) =  \sum_{k=1}^K P \left( f(X) = k, f(Y^I) = k\right),
$$
attains it maximum.

In this way, we look for the subset $I$ for which the original partition rule
applied to the less informative random vector $Y^I$  behaves as close as possible to
the procedure when is applied to the ``full information'' vector $X$. All components
with index outside the set $I$ are blinded in the sense that are constant.

In practice, the empirical version consist on the application of the next steps:

\begin{enumerate}
\item Given iid data $X_1, \ldots, X_n \in R^p$,   apply
the partition procedure to the data set and  obtain the empirical cluster allocation
function,
$$
f_n: R^p \rightarrow \{1, \ldots, K \},
$$
where now $f_n(x)$ is data dependent. The associated space partition will be denoted
by $G_k^{(n)} = f_n^{-1} (k)$, for $k = 1, \ldots, K$.

\item For a fixed value  $d < p$, given a subset of indices
$I \subset \{1, \ldots, p\}$, with $\#I = d$,  define the random vectors $ \{X^*_j,
1 \leq j \leq n \}$ verifying
$$
X^*_j [i] = X_j[i] \mbox{ if } i \in I, \mbox{ and } X^*_j[i] = \bar{X}[i]  \mbox{
otherwise},
$$
where $X[i]$ stands for the $i$--coordinate of the vector $X$, and $\bar{X}[i]$
stands for the $i$--coordinate of the average vector.

If we have used instead of the expected value other location parameter in the
population version (like the median), we substitute the average by the empirical
version  (the sample median).

\item  Calculate the {\em empirical objective function}
$$
h_n(I) = \frac{1}{n} \sum_{k=1}^K \sum_{j=1}^n \mathcal{I}_{\{f_n(X_j) =
k\}}\mathcal{I}_{\{f_n(X^*_j) = k\}},
$$
where $\mathcal{I}_A$ stands for the indicator function of the set $A$.

\item  Look for a subset $I_{d,n}=:I_n$, with
$\#I_n = d$,  that maximizes the empirical objective function $h_n$.
\end{enumerate}

\subsection{Consistency. Assumptions and main result}

As expected, the consistency of our variable selection procedure is linked to the
properties of the cluster partition method. We now give some conditions under which
our procedure is consistent.

\begin{description}
\item{{\em Assumption 1:}}

\noindent a) The partition procedure is strongly consistent, i.e., given $\epsilon >
0$, there exists a set $A(\epsilon) \subset R^p$ with $P(X \in A(\epsilon)) > 1 -
\epsilon$, such that for all $r>0$
$$
\lim_{n \rightarrow \infty} \sup_{x \in C(\epsilon, r)} |\mathcal{I}_{\{f_n(x) =
k\}} - \mathcal{I}_{\{f(x)=k\}}|=0 \mbox{ a.s., for } k=1, \ldots, K,
$$
where $C(\epsilon, r) = A(\epsilon) \cap B(0,r)$ stands for the intersection of the
set $A(\epsilon)$ and the closed ball centered at zero of radius $\epsilon$,
$B(0,r)$.

\noindent b)
$$d(X, \partial G_k^{n}) - d(X, \partial G_k) \rightarrow 0 \hspace{2mm} \mbox{a.s.,
for} \hspace{2mm} k=1, \ldots, K,$$ \indent where $d(X,\partial G_k)$ stands for the
distance from $X$ to the frontier of $G_k$.
\item{{\em Assumption 2:}}
$$\lim_{\delta \to 0} P(d(Y, \partial G_k) < \delta) = 0 \mbox{, for  }
k=1,\ldots, K.$$
\end{description}

\begin{figure}[tbp]
\begin{center}
\includegraphics[width=3.5cm]{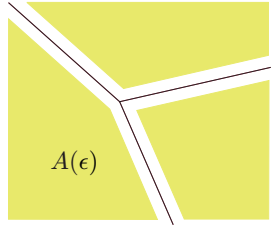}
\end{center}
\vspace{-15pt} \caption{Excluding a neighborhood of the partition boundaries, we
have almost sure uniform convergence of the function $f_n$ to $f$ over compact sets
(the color area  $A(\epsilon)$).}
\label{dibujito}
\end{figure}

Assumption 1a holds typically for cluster and classification rules, where the set
$A(\epsilon)$ is the complement of an $\epsilon$--neighborhood (``outer parallel
set") of the partition boundaries as shown in Figure~\ref{dibujito}, i.e.
$$
A(\epsilon)^c = \bigcup_{x\in \cup_{k=1}^K \partial G^k} B(x,\epsilon),
$$
where $B(x,\epsilon)$ denotes the ball with center $x$ and radius $\epsilon$.

\begin{theorem}
{ \rm\bfseries (Strong Consistency)} Let $ \{X_j: j \geq 1 \}$ be iid random vectors
with distribution $P_X$. Given $d$, $1 \leq d < p$, let $I_d$ be the family of all
subsets of $\{1, \ldots, p\}$ with cardinal $d$, and $\bf\it{I_{d,0}} \subset
\bf\it{I_d}$ the family of subsets where the maximum of $h(I)$ is attained,  for $I
\in \bf\it{I_d}$. Then, under assumptions 1 and 2 we have that there exists $n_0 =
n_0(\omega)$, such that
$$
I_n  \in \bf\it I_{d,0} \mbox{ for }  n \geq n_0(\omega) \mbox{ a.s.}
$$
\end{theorem}

The proof  is given in the Appendix.

\subsection{Selection of variables in simulated data}

In order to analyze our method performance, we carry out a Monte Carlo study for
some simulated date sets. In all of them we generated 100 observations in a three
dimensional variable space. The underlying distributions are mixtures of three
multivariate normals,

\begin{equation*}
\mbox{\boldmath $X$}={\footnotesize \left(
\begin{array}{c}
X_{1} \\ X_{2} \\ X_{3}
\end{array}
\right) }\sim \sum_{i=1}^{3}\alpha _{i}\,\mathcal{N}_{3}\left( \mbox{\boldmath
$\mu$}_{i},\mbox{\boldmath $\Sigma$}_{i}\right) ,
\end{equation*}
where $\alpha _{1}=\alpha _{2}=0.35$ and $\alpha _{3}=0.30$. The cluster structure
is defined through $X_{1}$ and $X_{2}$ and, to simplify,  we consider they are
independent in all the cases, with distributions given by
\begin{eqnarray*}
&X_{1}\sim \alpha _{1}\mathcal{N}(0,0.2)+\alpha _{2}\mathcal{N} (0.1,0.2)+\alpha
_{3}\mathcal{N}(0.9,0.2)&  \\ &X_{2}\sim \alpha _{1}\mathcal{N}(0,0.2)+\alpha
_{2}\mathcal{N} (0.9,0.2)+\alpha _{3}\mathcal{N}(0.1,0.2).&
\end{eqnarray*}
For the distribution of $X_3$ we consider two different scenarios.

\begin{figure}[tbp]
\begin{center}
\includegraphics[width=8cm]{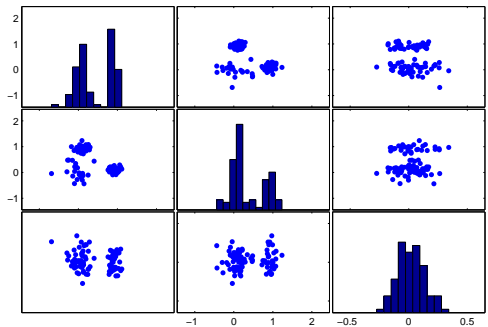}
\end{center}
\vspace{-15pt} \caption{Scatter plots and histograms from a three dimensional data
set generated following Case I description with $\sigma=0.2$} \label{3clusterplots}
\end{figure}

\noindent \emph{Case I:} $X_{3}$ is an independent ``noise'' variable with
distribution given by
$$
X_{3}\sim \mathcal{N}(0,\sigma),
$$
where $\sigma$ takes different values, $0.1$, $0.2$ and $0.3$.
Figure~\ref{3clusterplots} shows a simulated data set from these distributions with
$\sigma=0.2$. The three clusters are perfectly distinguish when plotting the pairs
$(x_{1},x_{2})$, however only two clusters are appreciated in the scatter plots that
consider $X_{3}$, as it is the case of the $X_{1}$ and $ X_{2}$ histograms.

\noindent \emph{Case II:} $X_{3}$ is not an independent variable and is given by
$$
X_{3}=(X_{1}+X_{2})/\sqrt{2}.
$$

In Table~\ref{tablasimul} we report the proportion of times where the information in
only one, two or three variables is enough to explain all the cluster allocations.
We also consider the effect of a possible reduction in the efficiency to only 95\%,
or 90\%, of correct allocations. In all the cases we carried out 1,000 replications
and follow the next steps.
\begin{enumerate}
\item Generate $\mbox{\boldmath $X$}_{1},\ldots ,\mbox{\boldmath $X$}_{100}$
observations.
\item Split the data into three cluster using the $k$--means algorithm.
\item Search the optimal subset of variables for 100\%, 95\% and 90\% efficiencies.
\end{enumerate}

\begin{table}
\centering
{\small \begin{tabular}{|c|cc|ccc|}
\hline &  &  & \multicolumn{3}{c|}{Number of variables } \\ &
&\multicolumn{1}{c|}{Efficiency} & 1 & 2 & 3 \\ \hline &  & 100\% & 0 & 0.997 &
0.003
\\ & $\sigma =0.1$ & 95\% & 0.005 & 0.995 & 0 \\ & &
90\% & 0.008 & 0.992 & 0 \\  \cline{2-6} &  & 100\% & 0 & 0.926 & 0.074 \\ $
X_{3}\sim \mathcal{N}(0,\sigma)$ & $\sigma =0.2$ & 95\% & 0.003 & 0.986 & 0.011 \\ &
& 90\% & 0.005 & 0.994 & 0.001
\\ \cline{2-6} &  & 100\% & 0 & 0.736 & 0.264 \\ & $\sigma=0.3$ &
95\% & 0 & 0.976 & 0.024 \\ &  & 90\% & 0.006 & 0.988 & 0.006 \\ \hline &  & 100\% &
0 & 0.146 & 0.854 \\ $ X_{3}={\displaystyle (X_{1}+X_{2})/\sqrt{2}} $  & & 95\% &
0.001 & 0.970 & 0.029
\\ &  & 90\% & 0.003 & 0.990 & 0.007 \\ \hline
\end{tabular}}
\caption{Simulation results from the Monte Carlo study carried out using the
distributions proposed in cases I and II} \label{tablasimul}
\end{table}

In the first case our variable selection method is very successful and selects only
the two variables $X_1$ and $X_2$ in almost all the simulations, for 100\%, 95\% and
90\% efficiencies. A different scenario appears with case II, where the third
variable is a linear combination of the first two variables. Only in the 14.6\% of
the times the two variable subset explains all the cluster allocations. This changes
dramatically when we allow for a 5\% or 10\% of miss--classified observations, now
97\% of the times the method selects only two variables, instead of three.

\begin{figure}[tbp]
\begin{center}
\includegraphics[height=3.5cm,width=4cm]{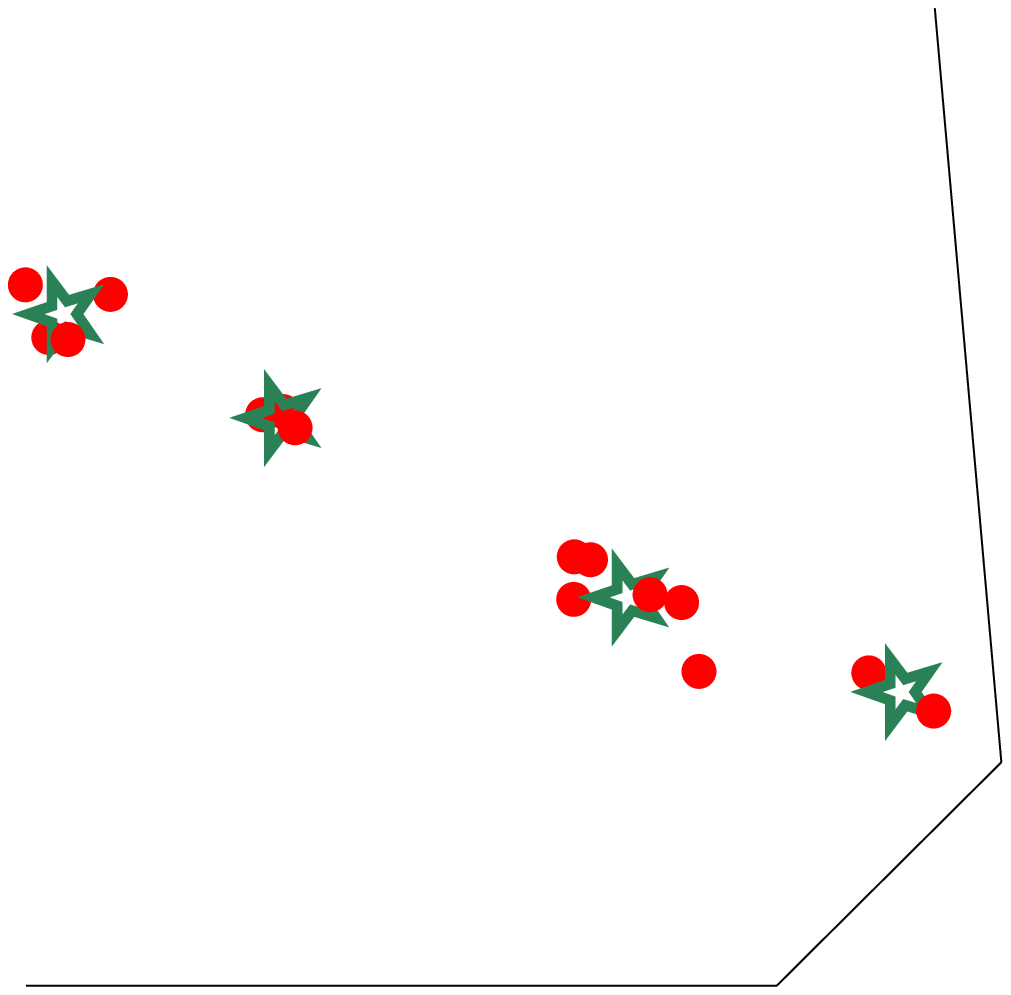} \hspace{1.5cm}
\includegraphics[height=3.5cm,width=4cm]{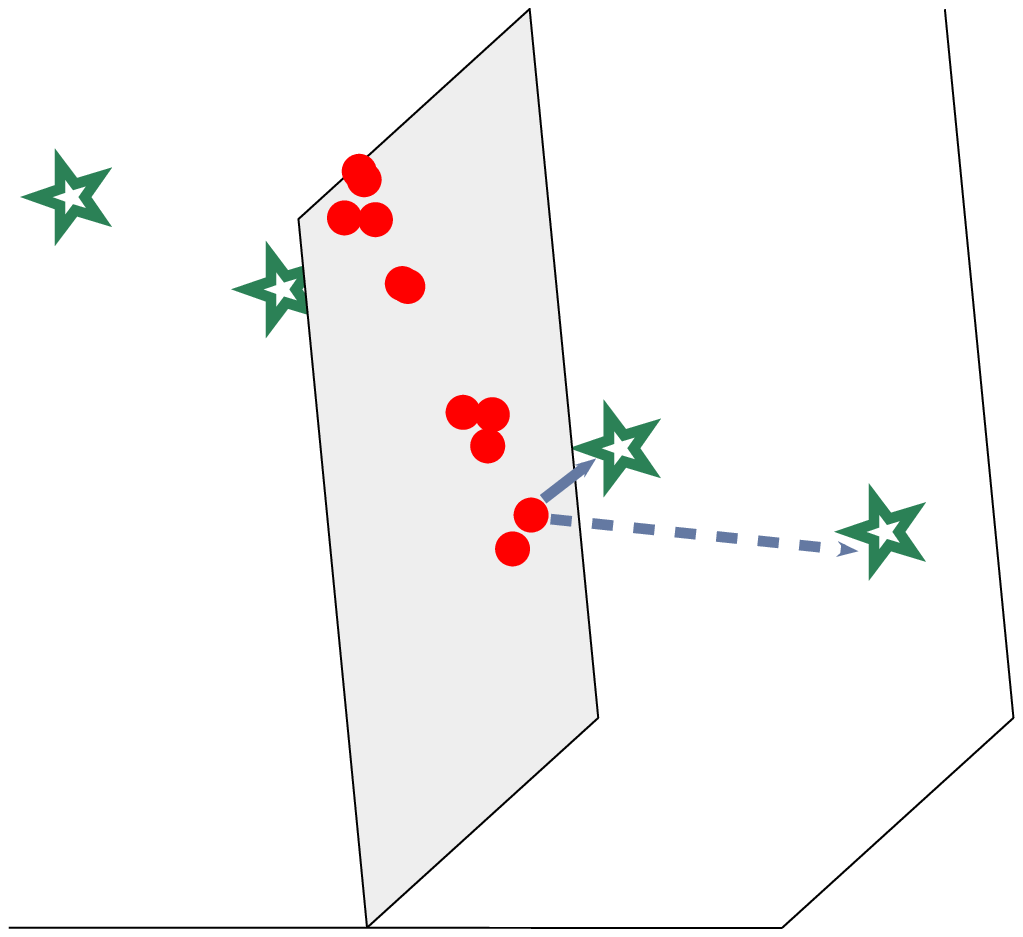}
\end{center}
\vspace{-15pt} \caption{a) Dots are the simulated data TSV05 as in Tadesse {\em et
al.} (2005) and  stars are the four $k$-means centers, b)  results of blinding the
vertical coordinate with the mean value.} \label{TSV05}
\end{figure}

Case II shows and interesting feature of the selection variable procedure, it is
able to eliminate noise variables, but it is unable to detect redundant information
from co--linear variables. This effect can be more clearly seen with the simulated
example proposed by Tadesse, Sha and Vannucci (2005). The data consists on the 15
three-dimensional observations displayed in Figure~\ref{TSV05}a. The first four
observations come from independent normals with mean $ \mu _{1}=5$ and variance
$\sigma _{1}^{2}=1.5$. The next three data come from independent normals with  mean
$\mu _{2}=2$ and variance $\sigma _{2}^{2}=0.1$. The following six data come from
independent normals with mean $\mu _{3}=-3$ and variance $\sigma _{3}^{2}=0.5$,
while the last two come from independent normals with mean $\mu _{4}=-6$ and
variance $\sigma _{4}^{2}=2$. Despite in Tadesse {\em et al.} (2005) the data set
was generated with twenty--dimensional observations instead of three--dimensional,
we call TSV05 to this data set.

We first run the $k$-mean algorithm with $k=4$, which classified correctly the whole
data set. Then, we run the variable selection procedure based on the mean value
(dropping out noisy non--informative variables). A closer look to this data
generating mechanism indicates that one should expect  to attain a 100\% efficiency
with only one variable, since we have the same cluster structure at the three
coordinates. However, the procedure was unable to find the cluster structure
blinding all variables except one. The efficiencies in Table~\ref{tablaTSV05a} show
that only the subset with the three variables classify all the data in their
original clusters. This result is expected since all the variables contain
information about the cluster, they are not noisy variables. However, as in case II,
these colinear variables are redundant and would be interesting to develop a
variable selection method able to detect them.

Figure~\ref{TSV05}b helps us to understand the main problem that appears when we
blind one variable by substituting all the data for the mean value. We observe the
case of blinding the vertical coordinate, that means a projection of all the data in
the shadow mean plane. As the mean is not a representative value for data generated
from a cluster structure, the allocations will be by chance to any of the clusters.
For instance, we point out the correct center for one projected data with a
discontinuous arrow, however in this case the closer center is a different one. This
data is wrongly allocated with the variable selection method. Remember that we blind
the variable but not the corresponding coordinate of the $k$-mean centers. Then to
eliminate not only noisy variables but also colinear variables the idea is to blind
with local information, instead of using the mean. This would not be a problem for
noisy variables and we will see in the next section that it is crucial for
multicolinearity.

\begin{table}[tbp]
 \centering {\small \begin{tabular}{|l||c|c|c|c|c|c|c|} \hline
Subset & $X_1$ & $X_2$ & $X_3$ & $X_1,X_2$ & $X_1,X_3$ & $X_2,X_3$ & $X_1,X_2,X_3$
\\ Efficiency & 60\% & 60\% & 60\% & 66.66\% & 73.33\% & 86.66\% & 100\%  \\ \hline
\end{tabular}}
\caption{Percentage of correct allocations in the TSV05 data set using the variable
selection method based on the mean.} \label{tablaTSV05a}
\end{table}

\section{Dealing with multicolinearity}

The previous procedure is mainly designed to find ``noisy'' non-informative
variables,  however as the simulated data set highlighted, it may fail in the
presence of colinearity.  In order to deal with this problem, we consider a quite
natural extension, changing the definition of the ``less informative'' vector $Y^I$.
Recall that we defined $Y^I_i = X_i$, if $ i \in I$, and $Y^I_i = E(X_i)$ otherwise.
Thus, for indices in the complement of the set $I$, $Y^I_i$ is defined as the best
constant predictor. Now the idea appears clearly,  to change means by conditional
means. We define the less informative vector $Z^I$ for indices $i$ in the complement
of the set $I$ as the conditional expectation of $X_i$ given the set of variables
$\{X_l: l \in I \}$, i.e. the best predictor of $X_i$ based on those variables. This
procedure will be able to deal with both kinds of problems. However, at a first
look, a shortcoming is that it will require a large sample size in order to estimate
the conditional expectation. Also the computational effort is quite bigger. The
choice of the smoothing parameter is also challenging, since it must involve not
more data than the size of the smaller cluster (if we think for instance in local
averages). If $m_n $  is the size of the smallest group for the partition procedure,
and for each $d$ and $n$, $r = r(n,d)$ is the number of nearest neighbor's we will
need to require that $r < m$, together with the standard conditions $r/n \rightarrow
0$  , and $n(r/n)^d \rightarrow \infty$, as $ n \rightarrow \infty$. We now describe
briefly the proposal in a precise setup.

\subsection{Population and empirical objetive function}

Given a subset of indices
$$
I = \{i_1 \ldots, i_d\} \subset \{1,\ldots, p\},
$$
let
$$
X[I] =: (X_{i_1}, \ldots, X_{i_d}), \hspace{2mm} \mbox{for} \hspace{2mm} i_1 < i_2 <
\ldots < i_d.
$$

We define the vector $Z^I:= Z = (Z_1, \ldots, Z_p)$, where $Z_i = X_i$ if $ i \in I$
and $Z_i = E(X_i| X[I])$ otherwise. Instead of the conditional expectation $E(X_i|
X[I])$ in order to attain robustness we can use local medians, or local M-estimates
(see for instance, Stone, 1977, Truong, 1989 or Boente and Fraiman, 1995).

For a fix integer $d < p$, now the {\em population objective function} is the set $I
\subset \{1, \ldots, p \}$, $\#I= d$,  for which the function
$$
h(I) =  \sum_{k=1}^K P \left( f(X) = k, f(Z^I) = k\right),
$$
attains it maximum.

In practice, the empirical version consists on the same steps than in the method
based on using the mean, except the second, that is substitute by the next step:

\begin{enumerate}
\item[2'.] For a fixed value of $d < p$, given a subset of indices $I \subset \{1, \ldots, p \}$, with $\#I = d$,
fix an integer value $r$ (the number of nearest neighbor to be used). For each $ j =
1, \ldots, n$, find the set of indices $C_j$ of the $r$-nearest neighbor's of
$X_j[I]$ among $\{ X_1[I], \ldots, X_n[I] \}$.

Now define the random vectors $ \{X^*_j, 1 \leq j \leq n \}$ verifying
$$
X^*_j [i] = X_j[i] \hspace{2mm} \mbox{if} \hspace{2mm} i \in I, \hspace{2mm}
\mbox{and} \hspace{2mm} X^*_j[i] = \frac{1}{r} \sum _{m \in C_j} X_m [i]
\hspace{2mm} \mbox{otherwise},
$$
where $X[i]$ stands for the $i$--coordinate of the vector $X$.

\end{enumerate}

A resistant procedure would take the local median instead of the local mean for $i
\notin I$ , i.e. $X^*_j[i] = median(\{X_m [i]: m \in C_j\})$.

\subsection{Consistency. Assumptions and main result}

As we have seen before the consistency of the variable selection method relays on
the properties of the cluster partition methods. Moreover, regularity conditions on
the boundary of the partitioning sets, in order to carry out the nonparametric
regression, are requested. Now, we give the conditions under which our procedure is
consistent. Together with Assumption 1 we will need:

\begin{description}
\item{{\em Assumption 3:}}
$$\lim_{\delta \rightarrow 0} P(d(Z, \partial G_k) < \delta) = 0, \mbox{ for all }
k=1,\ldots, K.$$

\item{{\em Assumption 4:}}
$$
\sup_x |g_{i,n}(x) - g_i (x) | \rightarrow 0 \hspace{2mm} a.s., \mbox{ for i }
\notin I
$$
where $g_i(x) = E(X_i | X(I) = x)$ is the corresponding non--parametric regression
functions and $g_{i,n}(x)$ is a consistent estimate of $g_i(x)$.
\end{description}

Assumption 4 allows to use any uniformly consistent estimate of the regression
function, although we have only describe above the case of $r$--nearest neighbor
estimates.

\begin{theorem}
{\rm\bfseries (Strong Consistency)} Let $ \{X_j: j \geq 1 \}$
be iid random vectors with distribution $P_X$. Given $d$, $1\leq d<p$, let $%
I_{d}$ be the family of all subsets of $\{1,\ldots ,p\}$ with cardinal $d$, and
$\mathit{{I_{d,0}}\subset {I_{d}}}$ the family of subsets where the maximum of
$h(I)$ is attained, for $I\in \mathit{{I_{d}}}$. Then, under assumptions 1, 3 and 4
we have that there exists $n_{0}=n_{0}(\omega )$, such that
$$
I_{n}\in \mathit{I_{d,0}\mbox{ for }n\geq n_{0}(\omega )\mbox{ a.s.}}
$$
\end{theorem}

The proof of Theorem 2 is very similar to that of Theorem 1 and we omit it in full
detail. We only point out the differences at the Appendix.

\subsection{TSV05 example revisited}

When we apply the conditional mean selection variable method to the 15 three
dimensional data of TSV05, we now obtain that with only one variable we attain a
100\% efficiency. This is the case of the second or the third variable, while with
the first one we obtain a 93.3\% efficiency.

A slightly different version of this example includes three new variables. The
additional noisy coordinates are generated from independent standard normal
distributions. The two variable selection procedures applied to the six--dimensional
data set produce exactly the same results as for the three dimensional data. The
noisy variables are not necessary to reach 100\% efficiency.


\section{A forward--backward algorithm}

A well known feature of the variable selection problem is the great number of
subsets that should be considered even for moderate values of $p$. An exhaustive
search guarantees to find the smaller subset of variables to achieve, at least, a
fixed percentage on the empirical objective function, however this procedure is non
feasible when many variables are considered. For instance, if $p=50$ we should check
among more than $10^{15}$ combinations. Alternatively, we propose a computationally
less expensive forward-backward search algorithm. We run the search meanly in the
forward mode and include the last step in the backward mode.

The algorithm starts from a one variable set and, progressively, includes new
variables with an iterative revision of the inclusions in each step. In general, the
backward search is less costly, but the leave-one-out strategy will make difficult
to find a small subset. When a set provides a percentage of good classifications
over the fixed percentage, the backward process starts the search of a more
parsimonious solution. To compute the objective function we can blind the variables
either replacing them by the mean or the conditional mean (in this case the
conditional distribution is towards the chosen subset up to that step). The
estimation of the conditional mean is done by nearest neighbors.

We distinguish three parts in the algorithm design, that are sequentially
implemented.

\begin{description}
\item[{\sf Part 1:}] Select the most ``influential'' variable $X^{\left( 1\right)}$
(the data assignation is more affected by its absence), blinding one by one all the
variables and selecting  the one with minimum value of the objective function,
\[
X^{\left( 1\right) }=\arg \min_{1\leq j\leq p}h_{n}\left( I_{j}\right) ,
\]
where $I_{j}$ denotes all the variables except the $j^{th}$, that is blinded.

\item[{\sf Part 2:}]
Sequential increment of variables one by one (forward search). In each step we look
for the accompanying variable, such that the new subset maximizes the number of
successfully data allocations. We also consider replacement of previously introduced
variables following the  iterative scheme described by Miller (1984) as a variation
of the classical forward-backward methods. The increment continues until the fixed
percentage of well classified data is reached.

\item[{\sf Part 3:}]
The subset is revised for unnecessary variables (backward search). The previously
introduced variables are questioned one by one whether they are necessary. The
algorithm stops when no further reduction can be found without loss of efficiency.

\end{description}

This procedure strongly depends on the order in the variable vector. To avoid (or
minimize) the label effect we run the algorithm for a random sample of the permuted
variables. We finally select the solutions that use the minimum number of variables.

In the following section we illustrate the algorithm performance in real data
examples. The matlab codes are available upon request to the authors.

\section{Real Data Examples}

\subsection{ Evaluation of educational programs}

We have survey data concerning education quality from 98 schools in the city and
suburbs of Buenos Aires (Argentine). The survey and posterior data analysis was
developed by Llach {\em et al.} (2006). An important objective on this study was to
find homogeneous groups of schools and the characterization of the clusters. The
selection variable method is a powerful tool to separate all the variables with real
influence from those that are non-informative.

At each school, a questionnaire with fifteen items  was fulfilled by the headmaster
and the teachers. The questions regard on the human and didactic resources, the
relationships between all the involved agents and the building physical condition.
All the answers range in a discrete scale between 1 and 100.  The items $V_1$ to
$V_8$ are answered by the headmasters and refer to their experience, aptitude,
school general knowledge, evaluation of the building conservation, evaluation of the
didactic resources, relationships with teachers, parents and students. The items
$V_9$ to $V_{15}$ are answered by the teachers and the questions are the same $V_1$
to $V_8$, except $V_3$ (school general knowledge) that is only answered by the
headmasters.

In Llach {\em et al.} (2006), a $k$--means cluster procedure was performed using the
98 fifteen dimensional vectors. The data were split into three clusters of sizes 45,
21 and 19 respectively.

The relationship between the clusters and the mean scores in a general knowledge
exam (GKE) and the mean socioeconomic level of the students (SEL) are shown in Table
\ref{Tabla GKE}. Both the GKE results and SEL are significatively different among
clusters, with ANOVA $p- \mbox{values}\ll 0.0001$. The clusters with higher mean
level of student knowledge correspond with those with also higher mean socioeconomic
level. The question now is which variables have relevant information to establish
this school grouping.

\begin{table}
\begin{center}
\begin{tabular}{|c||cc||cc|}
\hline & \multicolumn{2}{c||}{GKE} & \multicolumn{2}{c|}{SEL} \\
  & Mean & Std & Mean & Std \\
\hline Cluster 1 &  49.25 & 9.18  & 49.51 & 9.93 \\
       Cluster 2 &  58.04 & 13.01 & 63.49 & 16.39 \\
       Cluster 3 &  64.60 & 10.94 & 68.80 & 11.79 \\ \hline
\end{tabular}
 \caption{Means and standard deviations for the student general knowledge
 exams (GKE)  and the student socioeconomic levels (SEL).}
\label{Tabla GKE}
\end{center}
\end{table}

We select the variables that determine the clusters according to the first proposal,
with an exhaustive search (it is possible because of the moderate dimension of the
data). The clusters are completely explained (100\% efficiency) by  $ V_3, V_4, V_7,
V_8, V_{11}, V_{12}, V_{14}, V_{15} $, that are headmaster's school general
knowledge,  evaluation of the building conservation, relationships with parents and
students, and teacher's evaluation of the building conservation, evaluation of the
didactic resources, relationships with parents and students.

As the requested efficiency decays, so does the number of variables, the subset
includes the variables $V_1, V_4, V_7, V_{11}, V_{12}, V_{14}$ for 98\% efficiency.
For 95\% efficiency several subsets of size six were found. To achieve 92\%
efficiency, we found two optimal subsets with only four variables $V_4, V_7, V_{11},
V_{14}$ or $V_4, V_7, V_{12}, V_{14}$, that are headmaster's evaluation of the
building conservation, relationships between headmasters and parents, and between
teachers and parents. The elective variables are teacher's evaluation of the
building conservation or teacher's evaluation of the didactic resources In all the
cases the subsets contain information from the headmasters and the teachers.

With the aim of studying the algorithm performance  we run it with 100 permutations
and the results were consistent with the exact procedure. We found almost all the
subsets that were found before.

To refine the previous results we apply the conditional procedure to detect
colinearity. When we apply either the exact procedure or the algorithm we found the
same subsets, variables $V_2,V_3,V_4,V_7,V_8,V_{11},V_{12},V_{14}$ or
$V_3,V_4,V_7,V_8,V_9,V_{11},V_{12},V_{14}$ reach 100\% efficiency. For 97\%
efficiency the variables founded are $V_3,V_4,V_7,V_{11},V_{14}$, and only three
variables $V_4,V_7,V_{14}$ are requested to explain 91\% of the cluster allocations.
These final variables include headmaster's evaluation of the building conservation,
relationship between headmasters and parents, and between teachers and parents.

All over the data analysis we observe the importance of the relationships with the
parents, both with the headmasters and teachers. When we look for non--noisy
variables we found that also relationships with students have relevant information
about the cluster origin. However, these variables are eliminated from the final
subset when we use the conditional mean, this means that the opinion about the
relationships with the students contains redundant information.

\subsection{Identifying types of electric power consumers with functional data}

\begin{figure}
        \begin{center}
        \includegraphics[width=12.5cm]{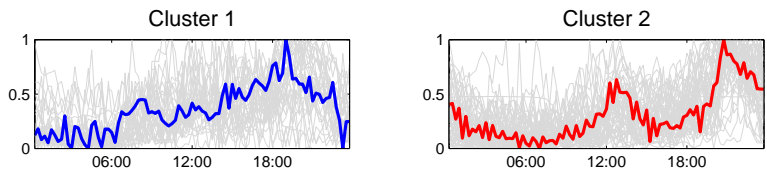}
        \caption{Home-consumers and cluster centers ($\alpha$-2-mean functions
        with trimming proportion $\alpha = 13/101$), of the two cluster structure for the
        electricity consume functional data.}
        \label{trim2means}
        \end{center}
\end{figure}

We consider the same example presented in Cuesta--Albertos and Fraiman (2006), where
an impartial--trimming cluster procedure is proposed for functional data. The study
was oriented to find  behavioral patterns of the electric power home--consumers in
the City of Buenos Aires. For each home, measurements were taken every 15 minutes
during all the weekdays of January 2001. The analyzed data were the vectors of
dimension 96 with the monthly averages  for a sample of 101 home--consumers. The
data were normalized in such a way that the maximum of each curve was equal to one.
Cuesta--Albertos and Fraiman (2006) found a two clusters structure, 13 outliers
apart. The resulting trimmed 2--mean functions (cluster centers) are shown in
Figure~\ref{trim2means}. Then the non-trimmed functions  were assigned to the
closest center and with this criteria the first cluster is composed of 33
home-consumers and the second one of 55. The remaining 13 data have been considered
as outliers.

In this example, the set of variables includes all the electricity consumptions in
the 15 minutes time--intervals in a day, that is 96 variables. A number too large
for the computation of the exact objective functions in all the possible subsets.
Therefore, in order to find the more relevant ``windows--times'' for the cluster
procedure we need to run the forward-backward search algorithm. We apply both the
mean and the conditional mean selection of variable algorithms for a 90\%, 95\% and
100\% of efficiency. For the calculation of the conditional mean, we consider 5, 10
and 33 nearest neighbor's (NN). The results after 100 permutations are summarized in
Table~\ref{Tablaluz}.

\begin{table}[b]
\begin{center}
{\small \begin{tabular}{|c||cc||cc||cc|} \hline
        &         & num. of   &        & num. of   &        & num. of   \\[-5pt]
 NN &  Effic. & variables & Effic. & variables & Effic. & variables  \\ \hline \hline
  5 &  90\%   &  4        & 95\%   &   6       & 100\%  &  9     \\
 10 &         &  3        &        &   5       &        & 16     \\
 33 &         &  3        &        &   7       &        & 28     \\ \hline \hline
 Mean &90\%   & 15        & 95\%   &  22       & 100\%  & 33     \\ \hline
\end{tabular}}
\caption{Optimal number of variables (time--intervals) for different
 efficiency percentages and number of nearest neighbor's using both the mean and
 the conditional mean selection of variable algorithms.} \label{Tablaluz}
\end{center}
\end{table}

\begin{figure}
        \begin{center}
        \includegraphics[width=8cm]{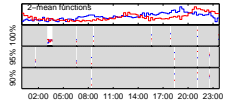}
        \caption{Two--mean electricity consume cluster centers for the functional data.
        Non--shadow time--intervals correspond to the subset of variables
        found by the 5-NN Conditional Mean Algorithm, with different
        degrees of efficiency.}
        \label{sol5NN}
        \end{center}
\end{figure}

The use of the conditional mean algorithm, instead of the faster mean algorithm,
reduces in all the cases the number of time--intervals that provides enough
information to characterize the two electric power home--consumer typologies. The
results show that the choice of the number of nearest neighbor's is also important,
although the method seems to be less sensitive than non--parametric regression.
However, it is an important problem to be solved. In our case, the results for 5-NN
are quite satisfactory: for a 100\% of efficiency, there is only one solution with 9
variables; for a 95\% of efficiency we found 15 different solutions with six
variables; while for a 90\% of efficiency we found 5 different solutions with four
variables. We choose one of them to illustrate in Figure~\ref{sol5NN} the
``window--times'' (non-shadow areas) which seems relevant.

\begin{figure}
        \begin{center}
        \includegraphics[width=13.5cm]{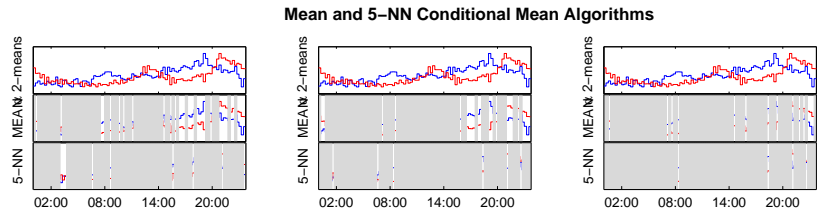}
        \caption{Two--mean electricity consume cluster centers for the functional data.
        Non--shadow time--intervals correspond to the subset of variables
        found by the Mean and the 5-NN Conditional Mean Algorithms, for different
        degrees of efficiency.}
        \label{meanand5NN}
        \end{center}
\end{figure}

The most informative consume registers are confined to a few ``window--times'' (see
Figure~\ref{meanand5NN}) and the two types of electric power consumers are mainly
characterized by their different aptitudes at some time--intervals in the morning
(7:00 to 11:00), evening (15:00 to 19:00), night (21:00 to 24:00) and early morning
(3:00 to 4:00). Comparing the mean and the 5-NN conditional solutions we observe
that the redundant information, specially at evening and night, is summarized in the
smaller subset of variables found by the mean conditional algorithm. When we reduce
the degree of efficiency and accept a number of missclassifications, the importance
of the early morning behavior diminished.

\section{Final Remarks}

We propose two variable selection procedures particularly design for partition rules
(typically supervised and un-supervised classification methods) that help to
understand the results for high-dimensional data. Both methods are strongly
consistent. The second procedure, based on conditional means, is much more flexible
and takes into account general dependence structures within the data. The
performance of our proposals in simulated and  real data examples is quite
impressive.

For low or moderate dimensional data an exhaustive search is possible for even the
case of 100\% efficiency. However, it is unfeasible for high-dimensional data and we
propose a forward-backward algorithm. We compare the algorithm performance with the
exhaustive search in some of the examples and the results are very positive since
they provide the same subsets. However, it will demand a considerable computational
effort in time that suggest that some additional research should be consider in this
aspect.

\section*{Acknowledgements}

Ricardo Fraiman and Ana Justel have been partially supported by the Spanish
Ministerio de Educaci\'on y Ciencia, grant MTM2004-00098.

\section*{Appendix}

\noindent {\em Proof of Theorem 1}

To simplify the proof, we assume that there exists a unique subset $I_{d,0} =: I_0 =
\{i_1, \ldots, i_d\} \subset \{1, \ldots, p\}$ that maximize $h(I)$ for $I \in
\bf\it{I_d}$. Then we point out the differences with the case of more than one
subset along the lines of the proof.

As the optimization is  over all the $d$ combinations of the $p$ variable indices, a
finite number, it suffices to show that
\begin{eqnarray}
\label{conv}  \lim_{n \rightarrow \infty} h_n(I) = h(I) \mbox{ a.s., for all } I \in
I_{d}.
\end{eqnarray}
Indeed, since there exist a unique set $I_0 \in \bf\it I_{p}$ that maximizes $h(I)$,
there exists $\eta>0$ such that
$$
h(I_0) > h(I) + \eta \hspace{2mm} \mbox{, for all } I \neq I_0, I \in I_{p}.
$$
We have from (\ref{conv}) that for all $I \in \bf\it I_{p}$
$$
|h_n(I) - h(I) | < \frac{\eta}{2}, \hspace{2mm} \mbox{if} \hspace{2mm} n \geq
n_0(I,\omega),
$$
which entails
\begin{eqnarray}
\label{conv1} h_n(I) < h(I) + \frac{\eta}{2} < h(I_0) - \frac{\eta}{2}, \hspace{2mm}
\mbox{if} \hspace{2mm} I \neq I_0.
\end{eqnarray}
Since we also have
$$
h_n(I_0) > h(I_0) - \frac{\eta}{2} > h(I) - \frac{\eta}{2},
$$
we conclude that $I_0$ maximizes $h_n(I)$ if $n \geq n_0(w)$ a.s.

If there exists more than one subset in $I_{d,0}$, the argument is the same by
replacing $I_0$ by $\bf\it I_{d,0}$.

Now, it remains to show (\ref{conv}), which reduces to prove that
\begin{eqnarray}
\label{ec1} \lim_{n \to \infty} \frac{1}{n} \sum_{j=1}^n
\mathcal{I}_{\{f_n(X_j)=k\}} \mathcal{I}_{\{f_n(X^*_j)=k\}} = P(f(X)=k, f(Y)=k)
 \hspace{2mm} \mbox{a.s.,}
\end{eqnarray}
for $k=1, \ldots, K.$

Finally, the equation (\ref{ec1}) will follow if we show that, for all fixed $k$,
\begin{eqnarray}
\label{ec11}
 \lim_{n \to \infty} \frac{1}{n} \sum_{j=1}^n \mathcal{I}_{\{f_n(X_j)=k\}}
 \mathcal{I}_{\{f_n(Y_j)=k\}} = P(f(X)=k, f(Y)=k)
 \mbox{ a.s.}
\end{eqnarray}
and
\begin{eqnarray}
\label{ec12}
 \lim_{n \to \infty} \frac{1}{n} \sum_{j=1}^n \mathcal{I}_{\{f_n(X_j)=k\}}
 [\mathcal{I}_{\{f_n(X^*_j)=k\}}- \mathcal{I}_{\{f_n(Y_j)=k\}}]
= P(f(X)=k, f(Y)=k) \mbox{ a.s.}
\end{eqnarray}

First we show that by the Assumption 1a we have
\begin{eqnarray}
\label{ec111} lim_{n \rightarrow \infty}  \frac{1}{n} \sum_{j=1}^n
{\mathcal{I}_{\{f_n(X_j)=k\}}}{\mathcal{I}_{\{f_n(Y_j)=k\}}} -
{\mathcal{I}_{\{f(X_j)=k\}}}{\mathcal{I}_{\{f(Y_j)=k\}}} = 0 \ \ \mbox{a.s.}.
\end{eqnarray}

The left hand side of (\ref{ec111}) is majorized by

$$
\frac{1}{n} \sum_{\{X_j \in C(\epsilon,r)\} \cap \{Y_j \in C(\epsilon,r)\}} |
{\mathcal{I}_{\{f_n(X_j)=k\}}}{\mathcal{I}_{\{f_n(Y_j)=k\}}} -
{\mathcal{I}_{\{f(X_j)=k\}}}{\mathcal{I}_{\{f(Y_j)=k\}}}| +
$$
$$
 \frac{1}{n} \sum_{\{X_j \notin C(\epsilon,r)\} \cup \{Y_j \notin C(\epsilon,r)\}} |
{\mathcal{I}_{\{f_n(X_j)=k\}}}{\mathcal{I}_{\{f_n(Y_j)=k\}}} -
{\mathcal{I}_{\{f(X_j)=k\}}}{\mathcal{I}_{\{f(Y_j)=k\}}}|
$$
The first term converges to zero for any $\epsilon$ and $r$ by Assumption 1a, while
the second term is dominated by
$$
\frac{1}{n} \# \{ 1 \leq j \leq n: \{X_j \notin C(\epsilon,r) \cup \{Y_j \notin
C(\epsilon,r)\} \}
$$
which converges a.s. to
$$
P(\left(\{X \notin C(\epsilon,r) \cup \{Y \notin C(\epsilon,r)\}\right).
$$
Since this last limit can be made arbitrarily small choosing $\epsilon$ and $r$
adequately, (\ref{ec111}) holds. Finally from the Law of Large Numbers we get
$$
\lim_{n \to \infty} \frac{1}{n} \sum_{j=1}^n
\mathcal{I}_{\{f(X_j)=k\}}\mathcal{I}_{\{f(Y_j)=k\}} = P(f(X)=k,f(Y)=k)  \mbox{
a.s.,}
$$
which concludes the proof of (\ref{ec11}).

For the proof of the equation (\ref{ec12}), the way in which the random vectors
$X^*_j$ and $Y_j$  have been defined, for a fixed subset $I$, implies that all the
$i$--coordinates of  $X^*_j - Y_j$ are zero for $i \in I$, while the rest of them
(for $i \notin I$) are given by
$$
\bar{X}[i] - E(X[i]).
$$
We recall that $X[i]$ stands for the $i$--coordinate of the vector $X$, and
$\bar{X}[i] = \frac{1}{n} \sum_{j=1}^n X_j[i]$.

The vectors $X^*_j - Y_j$ are all the same (i.e. the difference do not depend on
$j$), and are given by
\begin{eqnarray}
\label{ec121} (X^*_j - Y_j) [i] = (\bar{X}[i] - E(X[i])) \mathcal{I}_{\{i \notin
I\}},  \mbox{ for } j=1,\dots n.
\end{eqnarray}
>From the Law of Large Numbers we get
$$
\max_{ j = 1,\dots , n} ||X^*_j - Y_j|| \rightarrow 0 \hspace{2mm} \mbox{a.s.}
$$
The proof of (\ref{ec12}) will be complete if we show that
\begin{eqnarray}
\label{ec122} \# \left\{ j: f_n(X^*_j)=k, f_n(Y_j) \neq k, f_n(X_j) = k \right\} /n
\rightarrow 0 \hspace{2mm} \mbox{a.s.},
\end{eqnarray}
and
\begin{eqnarray}
\label{ec123} \# \left\{ j: f_n(Y_j)=k, f_n(X^*_j) \neq k, f_n(X_j) = k \right\} /n
\rightarrow 0 \hspace{2mm} \mbox{a.s.}
\end{eqnarray}
We  now define the  sets  B and C as follows:
$$
B = \left\{ \omega \in \Omega: \max_{ j = 1,\dots,n} ||X^*_j - Y_j|| \rightarrow 0
\right\},
$$

$$
C_j = \left\{\omega \in \Omega: d(X_j, \partial G^{(n)}_k) - d(X_j, \partial G_k)
\rightarrow 0 \right\}
$$
and
$$
C = \bigcap_{j=1}^{\infty} C_j.
$$
By the Assumption  1b we have that $P(B \cap C)=1$. Therefore, given $ \delta > 0$
and $\omega \in B \cap C$, there exists $n_0 = n_0(\omega, \delta)$ such that
$$
max_{j=1,\dots,n} ||X^*_j - Y_j|| \leq {\delta}/{2} .
$$
Given $\omega \in B \cap C$, we also have the following inclusions:
\begin{eqnarray*}
\left\{ j: f_n(X^*_j)=k, f_n(Y_j) \neq k, f_n(X_j) = k \right\} & \subset \left\{ j:
d(X^*_j, \partial G^{(n)}_k) < \delta \right\} \\ & \subset \left\{ j: d(Y_j,
\partial G_k) < 2\delta \right\},
\end{eqnarray*}
which imply that the left hand side of (\ref{ec122}) is majorized by
$$
\# \left\{ j: d(Y_j, \partial G_k) < 2\delta \right\}/n \leq \frac{1}{n}
\sum_{j=1}^n \mathcal{I}_{\{d(Y_j, \partial G_k) < 2\delta\}},
$$
which converges, as $n \rightarrow \infty$, to
$$
P(d(Y, \partial G_k) < 2\delta).
$$
Finally, from the Assumption 2   we get that
$$
\lim_{\delta \rightarrow 0 }P(d(Y, \partial G_k) < 2\delta)=0,$$ which concludes the
proof of (\ref{ec122}). The proof of (\ref{ec123}) is completely analogous.

\

\noindent {\em Proof of Theorem 2}

The proof goes on the same lines as the proof of Theorem 1. The only difference is
that now
$$
\max_{j=1,\dots,n} ||X^*_j - Z_j|| \rightarrow 0 \mbox{ a.s.}
$$
follows from Assumption 4.

\section*{References}

\begin{list}{}{\leftmargin.5cm\listparindent-.5cm} \item\hspace{-.2cm}
\vspace{-0.75cm}

Boente, G. and Fraiman, R. (1995), ``Asymptotic distribution of smoothers based on
local means and local medians under dependence''. \textit{Journal of Multivariate
Analysis}, \textbf{54}, 77--90.

Breiman, L., Friedman, J.H., Olsen, R.A. and Stone, C.J. (1984),
\textit{Classification and Regression Trees}. Wadsworth. Belmont, California.

Cuesta-Albertos, J. and Fraiman, R. (2006), ``Impartial trimmed k-means for
functional data''. \textit{Computational Statistics and Data Analysis} (in press).

Fowlkes, E.B., Gnanadesikan, R. and Kettenring, J.R. (1988), ``Variable selection
clustering''. \textit{Journal of Classification}, \textbf{5}, 205--228.

Fisher, R.A. (1936), ``The use of multiple measurements in taxonomic problems''.
\textit{Annals of Eugenics}, \textbf{7}, 179--188.

Fix, E. and J.L. Hodges (1951), ``Discriminatory analysis nonparametric
discrimination: Consistency properties''. Technical Report, U.S. Force School of
Aviation Medicine, Randolph Field, Texas.

Green, P.J. (1995), ``Reversible-Jump Markov Chain Monte Carlo computation and
Bayesian model determination''. \textit{Biometrika}, \textbf{82}, 711--732.

Hartigan, J.A. (1975), \textit{Clustering Algorithms}. John Wiley \& Sons, Inc. New
York.

Hern\'andez, A. and Velilla, S. (2005), ``Dimension Reduction in Nonparametric Kernel
Discriminant Analysis''. \textit{Journal of Computational and Graphical Statistics},
\textbf{14}, 847--866.

Llach, J.J. (2006). \textit{El desaf\'{\i}o de la equidad educativa: diagn\'ostico y
propuestas}. Ed. Granica, Buenos Aires.

Kaufman, L. and Rousseeuw, P.J. (1987), ``Clustering by means of Medoids''. In
\textit{Statistical Data Analysis Based on the $L_1$-Norm and Related Methods}, pp.
405--416. Y. Dodge. North-Holland.

MacQueen, J.B. (1967), ``Some Methods for classification and Analysis of
Multivariate Observations''. In \textit{Proceedings of 5-th Berkeley Symposium on
Mathematical Statistics and Probability}, Berkeley, University of California Press,
281--297.

Miller, A.J. (1984), ``Selection of subsets of regression variables''.
\textit{Journal of the Royal Statistical Society, A}, \textbf{147}, 389--425.

Pe\~na, D. and Prieto, F.J. (2001), ``Cluster Identification using Projections''.
\textit{Journal of American Statistical Association}, \textbf{96}, 1433--1445.

Stone, C. (1977), ``Consistent nonparametric regression'' (with discussion).
\textit{Annals of Statistics}, \textbf{5}, 595--645.

Tadesse, M.G., Sha, N. and Vannucci, M. (2005), ``Bayesian variable selection in
clustering high-dimensional data''. \textit{Journal of American Statistical
Association}, \textbf{100}, 602--617.

Truong, Y.K. (1989), ``Asymptotic properties of kernel estimators based on local
medians''. \textit{Annals of Statistics}, \textbf{17}, 606--617.

\end{list}

\end{document}